\documentclass[a4paper,10pt]{article}
\usepackage{amsfonts,latexsym,amssymb}

\usepackage[francais]{babel}

\newcommand{\la}{\langle\,}
\newcommand{\ra}{\,\rangle}

\renewcommand{\le}{\leqslant}
\renewcommand{\ge}{\geqslant}

\newtheorem{proposition}{Proposition}

\newtheorem{remark}{Remarque}

\begin{document}

\title{Sur deux propri\'et\'es de dualit\'e
des marches au hasard en
environnement al\'eatoire sur $\mathbb{Z}$}
\author{Didier Piau
\\
\begin{em}
Universit\'e Lyon 1
\end{em}
}
\date{}
\maketitle

\begin{abstract}
D'apr\`es Comets, Gantert et Zeitouni d'une part, et d'apr\`es
Derriennic d'autre part, certaines fonctionnelles associ\'ees \`a des
temps d'atteinte de marches au hasard en environnement al\'eatoire sur
${\mathbb Z}$ co\"\i ncident pour la marche elle-m\^eme et pour la
marche dans l'environnement renvers\'e.  Je montre que les deux
principes de dualit\'e ainsi exhib\'es sont alg\'ebriquement
\'equivalents, qu'ils d\'ecoulent de la seule propri\'et\'e de Markov
de la marche \`a environnement fix\'e et non pas de l'ergodicit\'e du
mod\`ele, 
et qu'on peut en donner des
versions finitistes et presque s\^ures.

\bigskip
\begin{center}
\begin{bf}
Abstract\end{bf}\end{center}

According to Comets, Gantert et Zeitouni on the one hand and to
Derriennic on the other hand, some functionals associated to the hitting
times of random walks in random environment on the integer line
coincide, for the walk itself and for the walk in the reversed
environment.
We show that these two duality principles are algebraically
equivalent, that they both stem from the Markov property of the walk
in a fixed environment, and not of the ergodicity of the model, and
that there exists finitist and almost sure versions of this duality.
\end{abstract}

\bigskip

\noindent
\textbf{English title:}
On two duality properties of random walks in random environment on the
integer line

\noindent
\textbf{Mots-cl\'es :}
Marches au hasard en environnement al\'eatoire, dualit\'es, fractions
continues.

\noindent
\textbf{Date :} Janvier 2000. R\'evision : juin 2001.

\noindent
\textbf{Classifications MSC :} 60K37 (60G50).

\section{Introduction}
\label{s.intro}

\paragraph{1. Cadre}
Soit $X$ une marche au hasard en environnement al\'eatoire sur
la droite $\mathbb{Z}$ des entiers.
Ainsi, $(E,{\cal E}, m)$ est
un espace de probabilit\'e, $S:E\to E$ une
bijection bime\-surable qui pr\'eserve $m$,
et $p:E\to(0,1)$ une application mesurable.
Pour $e\in E$ fix\'e, la suite $\{p_n(e)\,;\,n\in\mathbb{Z}\}$,
o\`u  $p_n:=p\circ S^n$,
permet de construire la loi $P_e$
d'une marche au hasard $X:=\{X_k\,;\,k\ge 0\}$ au plus proche voisin sur
$\mathbb{Z}$, dans l'environnement $e$, comme suit~: 
sous $P_e$, $X$ est une cha\^\i ne de Markov issue de
$X_0=0$ et, pour tout $k\ge 0$,
$$
P_e(X_{k+1}=n+1\,\vert\,X_k=n):=p_n(e)=:
1-P_e(X_{k+1}=n-1\,\vert\,X_k=n).
$$
Soit $q:=1-p$, $q_n:=1-p_n$, 
$E_e$ l'esp\'erance selon $P_e$ et
$\la\,\cdot\,\ra$ l'esp\'erance selon $m$.
Enfin, pour toute fonction $\phi$ d\'efinie sur $E$
et pour tout $n\in\mathbb{Z}$,
$\phi_n:=\phi\circ S^n$.

Les deux propri\'et\'es de dualit\'e que nous examinons concernent 
le retournement de $\mathbb{Z}$.
Pour des marches au hasard classiques
(donc, dans un environnement fix\'e), 
cette transformation donne des
relations parfois particuli\`erement simples. 
Par exemple, d'apr\`es H.~Dett, J.~Fill, J.~Pitman et 
W.~Studden~\cite{dfaps97}, 
voir aussi D.~Siegmund~\cite{sieg76},
si on intervertit les probabilit\'es des sauts $+1$ et $-1$ 
pour une marche au hasard sur la demi-droite $\mathbb{N}$ des entiers,
la fonction de Green $g_u$,
$u\in[0,1)$, de la
premi\`ere marche au hasard tu\'ee en une certaine barri\`ere et la fonction de
Green $h_u$ de la
deuxi\`eme marche au hasard 
r\'efl\'echie en cette m\^eme barri\`ere v\'erifient
$$
g_u(0,0)\cdot h_u(0,0)=(1-u)^{-1}.
$$
De m\^eme, d'apr\`es K.~Jansons~\cite{jans96}, si on inverse le sens de la
d\'erive d'une diffusion r\'eelle, les fonctions g\'en\'eratrices des
temps d'atteinte d'un m\^eme niveau par chacune des deux diffusions
obtenues sont reli\'ees par une \'equation simple.

\paragraph{2. Dualit\'e CGZ}
Dans leur preuve d'un principe de grandes d\'eviations portant sur
la vitesse $X_k/k$ de la marche au hasard en environnement al\'eatoire
d\'ecrite plus haut, et quand le d\'ecalage $S$ est ergodique, 
F.~Comets, N.~Gantert et O.~Zeitouni~\cite{cgz} 
montrent la relation de dualit\'e suivante.
Fi\-xons $u\in(0,1)$ et notons
$$
f(e):=E_e(u^{\tau_1}\,:\,\tau_1<+\infty),
\qquad
f'(e):=E_e(u^{\tau_{-1}}\,:\,\tau_{-1}<+\infty)
$$
les fonctions g\'en\'eratrices de $\tau_1$ et $\tau_{-1}$.
Ici, $\tau_n:=\inf\{k\ge 1\,;\,X_k=n\}$ est le premier temps de passage en
$n\in\mathbb{Z}$.
Alors, d'apr\`es la proposition~1 de la r\'ef.~\cite{cgz}, 
d\`es que $\log(p/q)$ est
$m$--int\'egrable,
\begin{equation}
\label{e.cgz}
\la\log f\ra=\la\log f'\ra+\la\log(p/q)\ra.
\end{equation}
Remarquons que la limite de (\ref{e.cgz}) quand $u\to1^-$
fournit l'\'egalit\'e
$$
\la\log P(\tau_1<+\infty)\ra=\la\log
P(\tau_{-1}<+\infty)\ra+\la\log(p/q)\ra.
$$
Par cons\'equent,
(\ref{e.cgz})
montre aussi que $\la\log g\ra=\la\log g'\ra$, avec
$$
g(e):=E_e(u^{\tau_1}| \tau_1<+\infty)
\quad\mbox{et}\quad
g'(e):=E_e(u^{\tau_{-1}}| \tau_{-1}<+\infty),
$$
puisque $g$ est le
rapport de $f$ \'evalu\'ee en $u$ et de $f$ \'evalu\'ee en $u=1$.
Si la marche au hasard est par exemple $m$--p.s.\ 
transiente vers $+\infty$, on en d\'eduit
$$
\la E(\tau_{-1}\,\vert\,\tau_{-1}<+\infty)\ra=\la E(\tau_{1})\ra,
$$
ainsi que l'\'egalit\'e des cumulants suivants.

\paragraph{3. Dualit\'e D}
Y.~Derriennic~\cite{derr},
en utilisant une repr\'esentation par cycles et poids,
r\'eussit \`a donner un crit\`ere de r\'ecurrence ou transience pour
des marches au hasard dont les sauts sont d'amplitude $+2$ ou $-1$ 
en environnement al\'eatoire ergodique sur $\mathbb{Z}$. 
La loi de $X$ est maintenant donn\'ee par
$$
P_e(X_{k+1}=n+2\,\vert\,X_k=n):=p_n(e)=:1-P_e(X_{k+1}=n-1\,\vert\,X_k=n).
$$
Un outil crucial est la relation de dualit\'e suivante,
\'enonc\'ee au lemme~3 de la r\'ef~\cite{derr}, 
dont nous changeons ci-dessous les notations.
Soit $c\ge0$ une fonction mesurable sur $E$.
Posons $r(0):=c=:r'(0)$ et, pour $n\ge 1$,
$$
r(n):=c_{-n}/r(n-1),
\qquad
r'(n):=c_{n}/r'(n-1).
$$
o\`u nous rappelons que $c_n:=c\circ S^n$.
Soient $x[n]$ et $x'[n]$ les deux fractions continues d\'efinies par 
$$
x[n]:=[r(0),r(1),\ldots,r(n)]
\quad\mbox{et}\quad
x'[n]:=[r'(0),r'(1),\ldots,r'(n)],
$$ 
et soient $x$
et $x'$ les fractions continues limites,
quand $n\to\infty$, de $x[n]$ et $x'[n]$.
Nous rappelons les notations qui concernent les fractions continues 
dans la section~\ref{s.equiv}.
\`A partir des propri\'et\'es de sym\'etrie des fractions
continues,
Derriennic montre que, d\`es que $\log c$ est $m$--int\'egrable,
\begin{equation}
\label{e.d}
\la\log x[n]\ra=\la\log x'[n]\ra
\end{equation}
pour tout $n\ge 0$ et donc, par convergence domin\'ee,
\begin{equation}
\label{e.dd}
\la\log x\ra=\la\log x'\ra.
\end{equation}

\paragraph{4. R\'esultats}
Nous montrons, sans l'hypoth\`ese d'ergodicit\'e, que l'extension suivante
de (\ref{e.cgz}) d\'ecoule de la seule propri\'et\'e de Markov de la
marche au hasard en environnement fix\'e.

\begin{proposition}
Soit $P_e$ la loi de la marche au hasard 
au plus proche voisin sur
$\mathbb{Z}$,
dans
l'environnement $e\in E$.
Pour $u\in(0,1)$ fix\'e et pour $n\ge 0$, soit
$$
f[n](e):=E_e(u^{\tau_1}\,:\,\tau_1<\tau_{-n-1}),
\qquad
f'[n](e):=E_e(u^{\tau_{-1}}\,:\,\tau_{-1}<\tau_{n+1}).
$$
Alors, pour tout $n\ge 1$,
\begin{equation}
\label{e.fexact}
f'[n]\cdot f[n-1]\circ S^n=f'[n-1]\cdot f[n]\circ S^n.
\end{equation}
Si $\log(p/q)$ est $m$--int\'egrable, $\log f[n]$ et $\log
f'[n]$ sont $m$--int\'egrables, et
\begin{equation}
\label{e.cgzext}
\la\log f[n]\ra=\la\log f'[n]\ra+\la\log(p/q)\ra.
\end{equation}
\end{proposition}

\begin{remark}
\begin{em}
La condition d'int\'egrabilit\'e de la proposition est \'equivalente
au fait que $\log (pq)$ est $m$--int\'egrable.  Remarquons aussi que 
$f[0]=p\,u$ et $f'[0]=q\,u$
donc (\ref{e.cgzext}) est \'evidente pour $n=0$.  Par ailleurs,
(\ref{e.cgz}) est la limite de (\ref{e.cgzext}) quand $n\to\infty$, 
donc (\ref{e.cgz}) se d\'eduit de (\ref{e.cgzext}) par
convergence monotone. On peut enfin souligner que (\ref{e.fexact})
n'entra\^\i ne pas de
relation $m$--presque s\^ure analogue pour les
limites $f$ et $f'$ puisque $f[n]$ et $f[n-1]$ apparaissent avec une
puissance du
d\'ecalage $S$ qui tend aussi vers l'infini.
\end{em}
\end{remark}

Nous montrons d'autre part que 
(\ref{e.d}) et (\ref{e.cgzext}) 
sont alg\'e\-bri\-que\-ment \'equi\-va\-len\-tes.  
Cette remarque nous
permet de d\'emontrer un analogue de (\ref{e.fexact}) pour $x[n]$ et
$x'[n]$, et donc de retrouver directement (\ref{e.d})--(\ref{e.dd}).


\paragraph{5. Plan}
La section \ref{s.markov} donne une preuve de (\ref{e.fexact}) et
(\ref{e.cgzext})
qui utilise uniquement la
propri\'et\'e de Markov.
Dans la section \ref{s.equiv}, nous montrons l'\'equivalence de 
(\ref{e.d}) et (\ref{e.cgzext}).
Enfin, la section \ref{s.ligne} montre comment une astuce donn\'ee
dans la r\'ef.~\cite{cgz} permet de
retrouver simplement (\ref{e.dd}),
mais sans doute pas (\ref{e.d}).

\paragraph{Remerciements}
Je tiens \`a remercier 
Nina Gantert et Yves Derriennic pour m'avoir fait conna\^\i tre
respectivement les r\'ef\'erences
\cite{cgz} et  \cite{derr} \`a l'origine de cette note,
et pour des discussions enrichissantes.

\section{Une preuve markovienne de (\ref{e.fexact})--(\ref{e.cgzext})}
\label{s.markov}
\noindent
Pour $n\ge 1$ fix\'e, la variable al\'eatoire $f'[n-1]\,f[n]\circ S^n$
ne d\'epend que des probabilit\'es de transition \`a partir de points
de $[0,n]$.  Notons donc $R$ la projection de $E$ sur $(0,1)^{n+1}$
d\'efinie par $R(e):=\{p_k(e)\,;\,0\le k\le n\}$.  

Alors,
$f'[n-1]\cdot f[n]\circ S^n=:F\circ R$ pour une certaine fonction
bor\'elienne $F$ (qui est m\^eme une fraction rationnelle).
On va montrer que $F$ est invariante par l'involution 
$I$ de $(0,1)^{n+1}$ qui envoie $\{p_k\,;\,0\le k\le
n\}$ sur $\{q_{n-k}\,;\,0\le k\le n\}$.
On peut lire $I$ comme le retournement du segment $[-1,n+1]$ qui
envoie $0$ sur $n$ et r\'eciproquement.
Puisque $F\circ I\circ R=f'[n]\cdot f[n-1]\circ S^n$, l'invariance par $I$  
entra\^\i nera bien (\ref{e.fexact}).

\`A cet effet,
notons $A$, $B$, $C$ et $D$ des fonctions bor\'eliennes telles que
\begin{eqnarray*}
A\circ R
& := & 
E^n(u^{\tau_{n+1}}\,:\,\tau_{n+1}<\tau_{-1})=f[n]\circ S^n,
\\
B\circ R 
& := & 
E^0(u^{\tau_{-1}}\,:\,\tau_{-1}<\tau_{n})=f'[n-1],
\\
C\circ R
& := &
E^0(u^{\tau_{n}}\,:\,\tau_{n}<\tau_{-1}),
\\
D\circ R
& := &
E^0(u^{\tau_{n+1}}\,:\,\tau_{n+1}<\tau_{-1}),
\end{eqnarray*}
o\`u $E^n$ d\'esigne 
l'esp\'erance pour une marche au hasard issue du point $n$.
Les fonctions compos\'ees $B\circ I$ et $C\circ I$ v\'erifient
\begin{eqnarray*}
B\circ I\circ R
& = & 
E^n(u^{\tau_{n+1}}\,:\,\tau_{n+1}<\tau_{0}),
\\
C\circ I\circ R
& = & 
E^n(u^{\tau_{0}}\,:\,\tau_{0}<\tau_{n+1}).
\end{eqnarray*}
La propri\'et\'e de Markov
au temps $\tau_0$ pour une marche au hasard issue de $n$ implique donc
$$
A=
B\circ I+C\circ I\cdot D.
$$
Il reste \`a remarquer que $D=C\cdot A$,
d'apr\`es la propri\'et\'e de Markov
au temps $\tau_n$ pour une marche au hasard issue de $0$,
pour
obtenir une expression de $A$ et donc de $F=A\cdot B$ en fonction de $B$
et $C$. 
Il vient
$$
F=B\cdot B\circ I/(1-C\cdot C\circ I),
$$
o\`u la division est licite puisque $C$ et $C\circ I\le u<1$.
Comme $B\cdot B\circ I$ et $C\cdot C\circ I$ sont invariants par
$I$,
$F$ l'est aussi.

Pour la preuve de (\ref{e.cgzext}), on remarque que l'int\'egrabilit\'e
de $\log p$ et $\log q$ entra\^\i ne celle de $\log f[0]$ et $\log
f'[0]$. 
Les suites $\{f[n]\,;\,n\ge 0\}$ et $\{f'[n]\,;\,n\ge 0\}$ sont
croissantes et major\'ees par $1$, donc tous les
$\log f[n]$ et $\log f'[n]$ sont int\'egrables.

La relation (\ref{e.fexact}) et l'invariance de $m$ sous l'action de
$S$ impliquent alors que $\la\log f[n]/f'[n]\ra$ ne d\'epend pas
de $n\ge 0$. Comme $f[0]/f'[0]=p/q$, on a d\'emontr\'e (\ref{e.cgzext}).


\section{L'\'equivalence de (\ref{e.d}) et (\ref{e.cgzext})}
\label{s.equiv}
\noindent
Soit $n\ge 1$,
$a:=p\,u=f[0]$ et $b:=q\,u=f'[0]$.
D'apr\`es la propri\'et\'e de Markov (simple),
$f[n]$ et $f'[n]$ v\'erifient
\begin{equation}
\label{ea.fg}
f[n]=a+b\,f[n]\,f[n-1]\circ S^{-1},\qquad
f'[n]=b+a\,f'[n]\,f'[n-1]\circ S^1.
\end{equation}
En effet, en conditionnant par le
premier pas de la marche au hasard, on a
par exemple
$$
f[n]=p\,u+q\,u\,E(u^{\tau_2}\,:\,\tau_2<\tau_{-n})\circ S^{-1},
$$
et il reste \`a calculer
$E(u^{\tau_2}\,:\,\tau_2<\tau_{-n})=f[n-1]\,f[n]\circ S^1$.

Avant d'\'etudier $x[n]$ et $x'[n]$, rappelons que, pour des r\'eels
$s_k$ donn\'es, tous positifs pour simplifier,
$[s_0,s_1,\ldots,s_n]$ d\'esigne une fraction continue, que l'on peut
d\'efinir par r\'ecurrence comme suit~: $[s_n]:=s_n$ et 
$$
[s_0,s_1,\ldots,s_n]:=s_0+1/[s_1,s_2,\ldots,s_n].
$$
On voit que, pour $\lambda\neq 0$,
$$
\lambda\,[s_0,s_1,s_2,s_3,\ldots]
=
[\lambda\,s_0,s_1/\lambda,\lambda\,s_2,s_3/\lambda,\ldots].
$$
Cette propri\'et\'e d'homog\'en\'eit\'e 
des fractions continues
est la seule que nous utilisons. Elle implique que
\begin{eqnarray*}
x[n]
& = & 
c\,[1,c_{-1},c_{-2}/c_{-1},\ldots]
\\
& = & 
c\,(1+1/[c,c_{-1}/c,c\,c_{-2}/c_{-1},\ldots]\circ S^{-1}).
\end{eqnarray*}
Ce calcul et un calcul analogue pour $x'[n]$ montrent que
\begin{equation}
\label{ea.xx}
x[n]=c\,(1+1/x[n-1]\circ S^{-1}),
\qquad
x'[n]=c\,(1+1/x'[n-1]\circ S).
\end{equation}
Bien s\^ur, les syst\`emes d'\'equations
(\ref{ea.fg}) et (\ref{ea.xx}) ne sont pas
identiques mais on passe de l'un \`a l'autre, au moins formellement,
en posant
$$
x[n]=-1/bf[n]\circ S^{-1},\quad
x'[n]=-1/a_{-1}f'[n],\quad
c=-1/a_{-1}b.
$$
Ou bien, dans l'autre sens, en posant par exemple
$$
f[n]=1/x[n]\circ S,\quad
f'[n]=-c/x'[n],\quad
a=1/c_1,\quad
b=-1.
$$
Cette transcription fait intervenir des quantit\'es
n\'egatives mais elle respecte, au contraire d'autres solutions
formellement correctes, 
les mesurabilit\'es de $x$, $x'$, $f$
et $f'$ (ainsi, $x$ est mesurable par rapport \`a $\{c_n\,;\,n\le
0\}$ et $f$  est mesurable par rapport \`a $\{p_n\,;\,n\le
0\}$), de m\^eme que les conditions initiales $f[0]=a$, $f'[0]=b$ 
et $x[0]=x'[0]=c$.
Ceci d\'emontre l'analogue suivant de (\ref{e.fexact}), valable
pour tout
$n\ge 1$~:
\begin{equation}
\label{e.xexact}
x'[n-1]\,x[n]\circ S^{n+1}=x'[n]\,x[n-1]\circ S^{n+1}.
\end{equation}
D\`es que tous les $\log x[k]$ sont $m$--int\'egrables,
l'invariance de la mesure sous l'action de $S$ assure donc que
$
\la\log x'[n]\ra= \la\log x[n]\ra,
$
pour tout $n\ge 0$, par r\'ecurrence sur $n$.
Or, $x[0]\le x[k]\le x[1]$ pour tout $k\ge 0$ et on peut v\'erifier que
$\log x[0]=\log c$ et 
$$
\log x[1]=\log c +\log(1+c_{-1})-\log c_{-1}
$$
sont tous deux $m$--int\'egrables si $\log c$ l'est.
On a d\'emontr\'e (\ref{e.d}).

\begin{remark}
\begin{em}
La m\'ethode utilis\'ee dans la r\'ef.~\cite{derr} 
permet \'egalement de retrouver
(\ref{e.xexact}). 
On peut pour cela utiliser la relation de sym\'etrie des fractions
continues rappel\'ee dans \cite{derr} et chercher une relation presque
s\^ure entre $x[n]$ et $x'[n]$, au lieu d'int\'egrer des \'egalit\'es
interm\'ediaires. 
On trouve alors pour $n\ge 0$,
$$
c\,x[n]\circ S^{n+2}\,\prod_{k=1}^{n}x[k]\circ S^{k+1}=c_2\,x'[n]\,
\prod_{k=1}^{n}x[k]\circ S^{k+2},
$$
relation qui entra\^\i ne bien (\ref{e.xexact}).
\end{em}
\end{remark}

\section{Les preuves ``en une ligne'' de (\ref{e.cgz}) et
(\ref{e.dd})}
\label{s.ligne}
\noindent
Nous rappelons la m\'ethode de la r\'ef.~\cite{cgz} pour d\'emontrer
(\ref{e.cgz}) et nous la copions pour montrer directement
(\ref{e.dd}).
Avec les notations de la Section \ref{s.equiv},
on remarque que $f$ et $f'$ v\'erifient les
relations de r\'ecurrence suivantes~:
$$
f=a+b\,f\,f\circ S^{-1},\qquad
f'=b+a\,f'\,f'\circ S.
$$
Plut\^ot que d'\'ecrire alors 
$f$ et $f'$ comme des fractions continues en
$\{p_n\,;\,n\le 0\}$ et $\{p_n\,;\,n\ge 0\}$ respectivement, 
on \'evalue simplement
\begin{eqnarray*}
b\,f\,(1-f\circ S^{-1}\,f') & = & 
b\,f-f'\,(f-a)=b\,f+a\,f'-f\,(b+a\,f'\,f'\circ S)
\\
& = & a\,f'\,(1-f\,f'\circ S).
\end{eqnarray*}
Si on pose $h=1-f\circ S^{-1}\,f'$, cette relation devient~: 
$b\,f\,h=a\,f'\,h\circ S$.
Ceci d\'emontre (\ref{e.cgz}) car $a/b=p/q$ et
$\la\log(h\circ S/h)\ra=0$.

\begin{remark}
\begin{em}
Puisque $f$ et $f'$ sont dans $[0,u]$,
$\log h$ est toujours $m$--int\'egrable.
La seule condition pour prouver
(\ref{e.cgz}) est donc bien que $\log(a/b)$ soit $m$--int\'egrable. 
Cette condition devient m\^eme inutile si on souhaite seulement
montrer que $<\log(bf/af')>=0$.
\end{em}
\end{remark}

Dans le m\^eme esprit que ci-dessus, 
une preuve de~(\ref{e.dd}) consiste \`a remarquer que $x$
et $x'$ v\'erifient les relations de r\'ecurrence~:
$$
x=c\,(1+1/x\circ S^{-1}),
\qquad
x'=c\,(1+1/x'\circ S).
$$
On introduit alors $y=c+x\,x'$ et 
on calcule
\begin{eqnarray*}
c_1\, x'\circ S\, y & = &
c_1\,x'\circ S\,c+c_1\,x\,c\,(1+x'\circ S)
\\
& = & 
c\,c_1\,
(x'\circ S+x)+c\,x'\circ S\,(x\,x\circ S-c_1)
\\
& = & c\,x\,(c_1+x'\circ S\,x\circ S)=c\,x\,y\circ S.
\end{eqnarray*}
Si $\log c$, $\log y$, $\log x$ et $\log x'$ sont $m$--int\'egrables,
on en d\'eduit
$$
\la\log(x'\circ S/x)\ra=\la\log(c/c_1)\ra+\la\log(y\circ S/y)\ra=0,
$$
ce
qui entra\^\i ne bien $\la\log x'\ra=\la\log x\ra$.
Or, en utilisant \`a nouveau les encadrements $x[0]\le x\le x[1]$ et
$x'[0]\le x'\le x'[1]$, on peut v\'erifier que 
$\log y$, $\log x$ et $\log x'$ sont bien $m$--int\'egrables
d\`es
que $\log c$ l'est.

\begin{remark}
\begin{em}
Cette m\'ethode ne semble pas redonner
les versions ``temps de sortie d'intervalles born\'es''
de (\ref{e.cgz}) et (\ref{e.dd}), 
c'est-\`a-dire (\ref{e.cgzext}) et (\ref{e.d}).
\end{em}
\end{remark}


\bigskip
\noindent
Universit\'e Claude Bernard Lyon 1
\\
LaPCS --
Laboratoire de Probabilit\'es, Combinatoire et Statistique
\\
43, boulevard du 11-Novembre-1918
\\
69622 Villeurbanne Cedex (France)
\\
[\smallskipamount]
{\tt Didier.Piau@univ-lyon1.fr}

\end{document}